\input amstex
\loadbold
\overfullrule0pt
\documentstyle{amsppt}
\magnification=1200
\hcorrection{.25in}
\advance\vsize-.75in
\topmatter
\title On the second moment for primes in an arithmetic progression\endtitle
\author D. A. Goldston$^{1}$ and C. Y. Y{\i}ld{\i}r{\i}m$^{2}$
\endauthor
\thanks $^1$ Research at MSRI is supported in part by
     NSF grant DMS-9701755, also supported by an NSF grant \endthanks
\thanks $^2$ Research at MSRI is supported in part by
     NSF grant DMS-9701755, also supported by a T\"{U}BITAK-NATO-B2 grant \endthanks
\address DAG:Department of Mathematics and Computer Science, San Jose
State University, San Jose, CA 95192, USA\endaddress
\email goldston\@mathcs.sjsu.edu \endemail
\address CYY:Department of Mathematics, Bilkent University, Ankara
06533, Turkey \endaddress
\email yalcin\@fen.bilkent.edu.tr \endemail
\leftheadtext{Goldston and Y{\i}ld{\i}r{\i}m}
\rightheadtext{Primes in Arithmetic Progressions}
\abstract
Assuming the Generalized Riemann Hypothesis, we obtain a lower bound within a constant factor of the conjectured asymptotic result for the second moment for primes in an individual arithmetic progression in short intervals. Previous results were averaged over all progression of a given modulus. The method uses a short divisor sum approximation for the von Mangoldt function, together with some new results for binary correlations of this divisor sum approximation in arithmetic progressions. 
\endabstract
\endtopmatter
\font\ger=eufm10
\def\gs{\hbox{\ger S}}
\def\a{\alpha}

\def\l{\lambda}
\def\L{\Lambda}

\def\de{\delta}

\def\g{\gamma}

\def\sc{\scriptstyle}

\def\d{\raise .5ex \hbox{\dag}}

\def\s2{\sum_{a=1}^q
  \bigl(\psi(x;q,a) \bigr)^2}
\def\star{\raise .5ex \hbox{*}}
\def\sumstar_#1{\setbox0=\hbox{$\scriptstyle{#1}$}
\setbox2=\hbox{$\displaystyle{\sum}$}
\setbox4=\hbox{${}\star\mathsurround=0pt$}
\dimen0=.5\wd0 \advance\dimen0 by-.5\wd2
\ifdim\dimen0>0pt
\ifdim\dimen0>\wd4 \kern\wd4 \else\kern\dimen0\fi\fi
\mathop{{\sum}\star}_{\kern-\wd4 #1}}
\def\sumprime_#1{\setbox0=\hbox{$\scriptstyle{#1}$}
\setbox2=\hbox{$\displaystyle{\sum}$}
\setbox4=\hbox{${}'\mathsurround=0pt$}
\dimen0=.5\wd0 \advance\dimen0 by-.5\wd2
\ifdim\dimen0>0pt
\ifdim\dimen0>\wd4 \kern\wd4 \else\kern\dimen0\fi\fi
\mathop{{\sum}'}_{\kern-\wd4 #1}}
\document
\TagsOnRight
\head 1. Introduction and Statement of results\endhead

In this paper we calculate a lower bound, of the same order
of
magnitude as conjectured, for the second moment
of primes in an arithmetic progression. Specifically
we examine
$$ I(x,h,q,a) := \int_x^{2x} \left(
\psi(y+h;q,a) -\psi(y;q,a) - {h\over \phi(q)}\right)^2\, dy
\tag 1.1 $$
where
$$\psi(x;q,a) = \sum_{\sc n\le x \atop \sc n\equiv a(\bmod\,
q)}
\Lambda(n), \tag 1.2$$
and $\Lambda$ is the von Mangoldt function.
We will take
$$ (a,q)=1, \,\,\quad x\geq 2, \,\,\quad 1 \leq q \leq h
\leq x, \tag 1.3 $$
(other ranges not being interesting).
We shall assume the truth of the Generalized Riemann
Hypothesis (GRH), which
implies, in particular,
$$ E(x;q,a) := \psi(x;q,a) - {x\over \phi(q)} \ll x^{{1\over
2}}\log^2 x
\,\, ,
\quad (q \leq x).  \tag 1.4$$
The idea of our method originates from the work of Goldston
\cite{3} for the
case of all primes, corresponding in the present formulation
to $q=1$. An
improved and generalized version of this result appeared in
\cite{5} as
\proclaim{Theorem A} Assume GRH. Then for any $\epsilon >0$
and
$1 \le {h\over q} \ll {x^{1\over 3}\over q^\epsilon \log^3
x}$ we have
$$\sum_{\sc a(\bmod\, q) \atop \sc (a,q)=1} I(x,h,q,a) \ge
{1\over 2
}xh\log
\Big(\big({q\over h}\big)^3x\Big) - O(xh(\log\log x)^3) .
\tag 1.5$$
Moreover, for almost all $q$ with $h^{3/4}\log^5 x\le q\le
h$ we have
$$I(x,h,q) \sim xh\log({xq\over h}). \tag 1.6$$
\endproclaim
For an individual arithmetic progression \"{O}zl\"{u}k
\cite{11} proved
unconditionally
\proclaim{Theorem B} For $1 \leq q \leq (\log
x)^{1-\delta}$, and
$h \leq (\log x)^c $ ($\delta$ and $c$ are any fixed
positive numbers)
satisfying  $q \leq h$, we have
$$ I(x,h,q,a) > ({1\over 2} - \epsilon){xh\over \phi(q)}\log
x  \tag 1.7 $$
for any $\epsilon$ and $x \geq X(\epsilon,c)$.
\endproclaim
We shall see below that the GRH implies a result of the type
in Theorem B
for much wider ranges of $q$ and $h$. An asymptotic estimate
for
$I(x,h,q,a)$
in certain ranges was shown by Y{\i}ld{\i}r{\i}m \cite{12}
to be implied by
GRH and a pair correlation conjecture for the zeros of
Dirichlet's
$L$-functions.
\proclaim{Theorem C} Assume GRH. Let $\a_1, \a_2, \eta$ be
fixed and
satisfying $ 0 < \eta < \a_1 \leq \a_2 \leq 1$, and let
$\delta = x^{-\a}$
where $\a_1 \leq \a \leq \a_2$. Assume, as $x \to \infty$,
uniformly for
$$ q \leq \min (x^{{1\over 2}}\delta^{{1\over 2}}\log^A x,
\delta^{-1}x^{-\eta}) \quad \quad
\hbox{($q$: prime or $1$)} \tag 1.8$$
and
$$ {x^{\a_{1}}\over \phi(q)}\log^{-3}x \leq T \leq
\phi(q)x^{\a_{2}}\log^{3}x \tag 1.9 $$
that,  for $(a,q) = 1$,
$$\sum_{\chi_1 , \chi_2 (\bmod \,
q)}\overline{\chi_{1}}(a)\chi_{2}(a)
\sum_{\sc 0 < \g_{1},\g_{2} \leq T \atop
{\sc L({1\over 2}+ i\g_{1}, \chi_{1}) = 0 \atop
\sc L({1\over 2}+ i\g_{2}, \chi_{2}) = 0 }}
x^{i(\g_1 - \g_2 )}{4\over 4 + (\g_1 - \g_2)^2} \,\,
\sim \, \phi(q){T\over 2\pi}\log qT .
\tag 1.10 $$
Then
$$\int_x^{2x} \left(
\psi(u+u\de ;q,a) -\psi(u;q,a) - {u\de\over
\phi(q)}\right)^2\, du \,\,
\sim {3\over 2}{\de x^{2}\over \phi(q)}\log {q\over \de}
\tag 1.11 $$
uniformly for $x^{-\a_{2}} \leq \de \leq x^{-\a_{1}} $ and
$q$ as in (1.8).
\endproclaim
It was also shown in \cite{12} that the left-hand side of
(1.10) is
$ \sim \phi(q) {T\over 2\pi}\log x $ for $ 1 \leq q \leq
x^{{1\over 2}}
\log^{-3}x $ when $ {x\over q}\log x \leq T \leq
e^{x^{{1\over 4}}}$.
These asymptotic values are what the diagonal terms
($\chi_1 = \chi_2 $) would contribute,
so the assumption (1.10) is a way of expressing that the
zeros of different Dirichlet $L$-functions are uncorrelated.
Theorem C is a generalization of
one half of a
result of Goldston and Montgomery \cite{4} for the case
$q=1$, where
an equivalence between the pair correlation conjecture for
$\zeta(s)$ and
the
second moment for primes was established. Since the argument
in \cite{4}
works
reversibly, a suitable converse to Theorem C is also
provable. The
restriction
to prime $q$ was made in order to avoid the presence of
imprimitive
characters.
The formula (1.11) involving differences $u\delta$ which
vary with $u$
can be converted to a formula involving a fixed-difference $h$.

Our main result is the following theorem.
\proclaim{Theorem 1} Assume GRH. Then for any $\epsilon > 0$
and
$$ q \leq h \leq (xq)^{{1\over 3}-\epsilon}, \tag 1.12 $$
we have
$$I(x,h,q,a) \geq {1\over 2}{xh\over \phi(q)}\log ({xq\over
h^3}) -
O({xh\over \phi(q)}(\log\log 3q)^3 ). \tag 1.13
$$
\endproclaim
\noindent Notice that the conditions in (1.12) imply that
both $h$ and $q$ are
$\ll x^{{1\over 2}-\epsilon}$.

The proof of the theorem uses some new results on the function
$\lambda_R(n)$ used as an approximation for the von Mangoldt
function in our earlier work. Propositions 2, 3, and 4 
embody these results, and we expect they will have further
applications to other problems.

\head 2. Preliminaries \endhead

We shall need the following in our calculations. Let
$$ f(n,x,h) = \int_{[x,2x]\cap [n-h,n)} 1 \, dy =
\cases n-x, &\text{for $x\le n<x+h$}\\
h,&\text{for $x+h\le n\le 2x $}\\
2x-n+h, & \text{for $2x<n\le 2x+h$}\\
0,& \text{elsewhere}. \endcases \tag 2.1 $$
\proclaim{ Lemma 1} For real numbers $a_n$ and $b_n$ we have

$$\align \int_x^{2x}\left( \sum_{y<n\le y+h}a_n\right)&
\left(  \sum_{y<m\le y+h}b_m  \right) \, dy =
\sum_{x<n\le 2x+h} a_n b_n f(n,x,h) \\
+ & \sum_{0<k\le h}\left(\sum_{x<n\le 2x+h-k}
(a_nb_{n+k} + a_{n+k}b_n) f(n,x,h-k)\right).
\endalign $$
\endproclaim
\proclaim{Lemma 2} Let $C(x) = \sum_{n\le x}c_n$. Then we
have
$$\sum_{x<n\le 2x+h}c_nf(n,x,h) =
\int_{2x}^{2x+h}C(u)\, du -
\int_{x}^{x+h}C(u)\, du .  $$
\endproclaim
Lemma 1 and Lemma 2 were proved in \cite{5}. We take this opportunity
to correct a minor error in Lemma 1 of \cite{5}.
In that lemma an extraneous term $h(c_{x+h}-c_{2x})$
was incorrectly included and should be removed. This
term then contributed an unnecessary error term in equations
(2.7),(2.14), and
(2.15) of [5]. However these same error terms correctly occurred for a
different reason in equation (2.9)
so that starting with equation (2.16) these error terms were correctly
included in the rest of [5].

Calling $\Delta\psi = \psi(y+h;q,a) - \psi(y;q,a) $ for
brevity, we
have from (1.1)
$$ I(x,h,q,a) = \int_{x}^{2x}(\Delta\psi)^2\, dy -
{2h\over \phi(q)}\int_{x}^{2x}(\Delta\psi)\, dy +
{h^2 x\over \phi^2 (q)}. $$
By the above lemmas and (1.4) we obtain
$$\align \int_{x}^{2x}(\Delta\psi)\, dy & =
\sum_{\sc x<n\le 2x+h \atop \sc n\equiv a (\bmod\, q)}
\Lambda(n) f(n,x,h)
\\
& = {xh\over \phi(q)} + \int_{2x}^{2x+h}E(y;q,a)\, dy -
\int_{x}^{x+h}E(y;q,a)\, dy \\
& = {xh\over \phi(q)} + O(hx^{{1\over 2}}\log^2 x) ,
\endalign $$
so that
$$ I(x,h,q,a) =  \int_{x}^{2x}(\Delta\psi)^2\, dy - {x
h^{2}\over \phi^2
(q)}
+ O({x^{{1\over 2}}h^{2}\log^2 x\over \phi(q)}). \tag 2.2$$
The integral $\int (\Delta\psi)^2$ leads to sums of the sort
$\sum \Lambda(n)\Lambda(n+k)$ which are in the territory of
the twin prime
conjecture.
In the  uninteresting case $1 \leq h \leq q$, only the sum
$\sum \Lambda^2 (n)$ is present, giving easily the
evaluation
$$ I(x,h,q,a) = {xh\over \phi(q)}\log x - {xh\over \phi(q)}
-
{x h^{2}\over \phi^2 (q)} + O(x^{{1\over 2}}h
\log^3 x), \, \quad (h \leq q \leq x). \tag 2.3$$

Now let $\lambda_{R}(n)$ be any arithmetical function, and
set
$$\psi_{R}(y;q,a) =
\sum_{\sc n \leq y \atop \sc n \equiv a (\bmod\,
q)}\lambda_{R}(n) \, \, ;
\quad \quad \Delta\psi_{R} = \psi_{R}(y+h;q,a)-
\psi_{R}(y;q,a). \tag 2.4$$
Trivially $\int (\Delta\psi - \Delta\psi_R )^2 \geq 0 $, so
that
$$ \int_{x}^{2x}(\Delta\psi)^2\, dy \geq
2\int_{x}^{2x}(\Delta\psi)\cdot (\Delta\psi_R )\, dy
- \int_{x}^{2x}(\Delta\psi_{R})^2\, dy . \tag 2.5$$
We apply Lemma 1 and Lemma 2 to these integrals to obtain
$$\align \int_{x}^{2x}(\Delta\psi)\cdot &(\Delta\psi_R )\,
dy \,\,  =
\sum_{\sc x < n \leq 2x+h \atop \sc n \equiv a (\bmod\, q)}
\Lambda(n)\lambda_{R}(n) f(n,x,h) \tag 2.6\,a \\
& + \sum_{\sc 0 < k \leq h \atop \sc q \mid k}
\sum_{\sc x < n \leq 2x+h-k \atop \sc n \equiv a (\bmod\,
q)}
[\Lambda(n)\lambda_{R}(n+k)+ \Lambda(n+k)\lambda_{R}(n)]
f(n,x,h-k) \\
= & \left(\int_{2x}^{2x+h} - \int_{x}^{x+h}\right)
\sum_{\sc n \leq u \atop \sc
n \equiv a (\bmod\, q)}\Lambda(n)\lambda_{R}(n)\, du \\
& + \int_{2x}^{2x+h} \sum_{0 < |j| \leq {u-2x\over q}}
\sum_{\sc N_{1} < n \leq N_{2} \atop \sc
n \equiv a (\bmod\, q)}\lambda_{R}(n)\Lambda(n+jq) \, du \\
& - \int_{x}^{x+h} \sum_{0 < |j| \leq {u-x\over q}}
\sum_{\sc N_{1} < n \leq N_{2} \atop \sc
n \equiv a (\bmod\, q)}\lambda_{R}(n)\Lambda(n+jq)\, du \, ,
\tag 2.6\,b
\endalign $$
where $N_1 = \max (0,-jq)$ and  $ N_2 = \min (u,u-jq)$.
Similarly
$$\align \int_{x}^{2x}(\Delta\psi_R )^2 \, dy \,\, = &
\sum_{\sc x < n \leq 2x+h \atop \sc n \equiv a (\bmod\, q)}
\lambda_{R}^2 (n) f(n,x,h)\\
&  + 2 \sum_{\sc 0 < k \leq h \atop \sc q \mid k}
\sum_{\sc x < n \leq 2x+h-k \atop \sc n \equiv a (\bmod\,
q)}
\lambda_{R}(n)\lambda_{R}(n+k) f(n,x,h-k) \\
= & \left(\int_{2x}^{2x+h} - \int_{x}^{x+h}\right)
\sum_{\sc n \leq u \atop \sc
n \equiv a (\bmod\, q)}\lambda_{R}^2 (n)\, du \\
& + 2 \int_{2x}^{2x+h} \sum_{0 < j \leq {u-2x\over q}}
\sum_{\sc n \leq u-jq \atop \sc
n \equiv a (\bmod\, q)}\lambda_{R}(n)\lambda_{R}(n+jq)\, du
\\
& - 2 \int_{x}^{x+h} \sum_{0 < j \leq {u-x\over q}}
\sum_{\sc n \leq u-jq \atop \sc
n \equiv a (\bmod\, q)}\lambda_{R}(n)\lambda_{R}(n+jq)\, du
.
\tag 2.7 \endalign$$

\head 3. The choice of $\l_{R}(n)$ and some number-theoretic
sums \endhead

As the auxiliary function we use
$$ \l_{R}(n) :=
\sum_{r \leq R}{\mu^2 (r)\over \phi(r)}
\sum_{d \mid (r,n)}d\mu(d) . \tag 3.1 $$
This function is known (\cite{3}, \cite{8}) to exhibit
behavior similar to
$\L(n)$ when considered on average in arithmetic
progressions, and it has
been
employed in related problems (\cite{1}, \cite{2}, \cite{5},
\cite{6}). An
upper
bound for $\l_{R}(n)$ is
$$|\l_{R}(n)| \leq
\sum_{d \mid n}d \sum_{\sc r \leq R \atop \sc d \mid
r}{1\over \phi(r)}
\leq \max_{r \leq R}({r\over \phi(r)})
\sum_{d \mid n}d \sum_{\sc r \leq R \atop \sc d \mid
r}{1\over r}
\ll d(n)\log R \log\log R . \tag 3.2 $$
To evaluate the sums which arise when (3.1) is used in
(2.6) and (2.7) we shall need some lemmas. In the following
$p$ will denote
a
prime number.
\proclaim{Lemma 3} {\rm (Hildebrand \cite{9})} We have for
each positive
integer $k$, uniformly in $R \geq 1$,
$$L_{k}(R) := \sum_{\sc n\leq R \atop \sc (n,k)=1}{\mu^2
(n)\over \phi(n)}
= {\phi(k)\over k}(\log R + c + v(k)) + O({w(k)\over
\sqrt{R}}),
\tag 3.3 $$
where
$$ \align c := & \gamma + \sum_{p}{\log p\over p(p-1)} ;
\quad
v(k):= \sum_{p\mid k}{\log p\over p}; \\
w(k):= & \sum_{d\mid k}{\mu^2(d)\over \sqrt{d}} =
\prod_{p\mid k}(1+{1\over \sqrt{p}}) ; \quad
v(1) = 0 ,\,  w(1) = 1.\tag3.4 \endalign $$
\endproclaim
\proclaim{Lemma 4} We have
$$  v(k) \ll \log\log 3k, \tag 3.5 $$
$$\sum_{p\mid k} {1\over \sqrt{p}} \ll  {\sqrt{\log k}\over
\log\log 3k},
\tag 3.6 $$
and
$$g(k) := \prod_{p\mid k}(1+{p\over p-1}) \ll 2^{\nu(k)}
(\log\log 3k)\tag 3.7 $$
\endproclaim
\demo{Proof} We show (3.7); the other inequalities can be
proved similarly.
Let $\nu(k)$ be the number of distinct prime factors of $k$,
which satisfies
the bound $\nu(k) \ll \displaystyle{\log k\over \log\log k}$.
We have
$$ \align
\log g(k) & = \sum_{p\mid k}\log (2+{1\over p-1})\\
& < \nu(k) \log 2 + \sum_{p\mid k}{1\over p} \\
& < \nu(k) \log 2 + \sum_{p \leq 2\log 2k}{1\over p} \\
& = \nu(k) \log 2 + \log\log\log 21k + O(1) ,
\endalign $$
where the prime number theorem and Mertens' theorem have
been employed. Exponentiating both sides we obtain (3.7).

\enddemo
\proclaim{Lemma 5} We have
$$ \sum_{\ell\mid k}{\mu^2 (\ell)\over \phi(\ell)}\log \ell
= {k\over \phi(k)}v(k).
\tag 3.8 $$
\endproclaim
\demo{Proof}
$$ \align
\sum_{\ell\mid k}& {\mu^{2}(\ell)\over \phi(\ell)}\log \ell
=
\log\prod_{\ell\mid k}
\ell^{\displaystyle{\mu^{2}(\ell)\over \phi(\ell)}} =
\log\prod_{p\mid k}p^{\,\,\,\displaystyle\sum_{\ell\mid k,
p\mid\ell}
{\mu^{2}(\ell)\over \phi(\ell)}} \\
& =
\log\prod_{p\mid k}p^{\displaystyle{1\over \phi(p)}
\displaystyle
\sum_{\ell'\mid {k\over p}, (\ell',p)=1}
{\mu^{2}(\ell')\over \phi(\ell')}}  =
\log\prod_{p\mid k}p^{\displaystyle{1\over \phi(p)}
\displaystyle{1\over 1+{1\over
\phi(p)}}\displaystyle\prod_{p'\mid
k}(1+\displaystyle{1\over \phi(p')})} \\
& = \log\prod_{p\mid k}p^{\,\,\displaystyle{1\over p}{k\over
\phi(k)}}
= {k\over \phi(k)}
\sum_{p\mid k}{\log p\over p} .
\endalign $$
\enddemo
\proclaim{Lemma 6} {\rm (Goldston \cite{3})} We have
$$\sum_{r \leq R}{\mu^2 (r)\over \phi(r)}
\sum_{\sc d \mid r \atop \sc (d,k)=1}{d\mu(d)\over \phi(d)}
=
\gs(k) + O({k\,d(k)\over R\phi(k)}), \tag 3.9 $$
where
$$ \gs(k) = \left\{ \eqalign{
       2C\prod_{{\sc p \vert k\atop \sc p>2}}\left({p-1\over
p-2}\right), &\quad
        \hbox{if $k$ is  even,
        $k\neq 0$;}\cr
0, \qquad  &\quad\hbox{if $k$ is  odd;}\cr}\right. \tag 3.10
$$
with
$$ C = \prod_{p>2}\left( 1 - {1\over (p-1)^2}\right). \tag
3.11 $$
\endproclaim
\demo{Proof} The proof can be found in \cite{3};
we just note that
$$ \sum_{\sc d\mid r \atop \sc (d,k)=1}{d\mu(d)\over
\phi(d)} =
{\mu(r)\over \phi(r)}\mu((r,k))\phi((r,k)) \, , \tag 3.12 $$
so the left-hand side of (3.9) may be expressed as
$$ \sum_{r=1}^{\infty}{\mu(r)\mu((r,k))\phi((r,k))\over
\phi^2 (r)}
+O(\sum_{r>R}{\mu^2 (r)\mu^2 ((r,k))\phi((r,k))\over \phi^2
(r)}). \tag
3.13$$
Here the first sum is $\gs(k)$ and the error term is $\ll$
the $O$-term in
(3.9).
\enddemo
\proclaim{Lemma 7} {\rm (Goldston and Friedlander
\cite{1})}We have
$$\sum_{0 < j \leq {h\over q}}(h-jq)\gs(jq) = {h^{2}\over
2\phi(q)} -
{h\over 2}\log {h\over q} + O(h(\log\log 3q)^3) . \tag 3.14
$$
\endproclaim
\proclaim{Lemma 8} {\rm (Hooley \cite{10})} Assuming GRH, we
have
$$ \sum_{\sc a (\bmod\, q) \atop \sc (a,q)=1} \max_{u \leq
x}|E(u;q,a)|^2
\ll x\log^4 x; \, \quad \hbox{for}\quad  q\leq x . \tag 3.15
$$
\endproclaim
\proclaim{Lemma 9}We have
$$ \align
\sum_{r \leq R}{\mu^2 (r)\sigma(r)\over \phi(r)} & \ll R
\tag 3.16 \\
\sum_{r \leq R}{\mu^2 (r)\sigma_{{1\over 2}}(r)\over
\phi(r)}
& \ll \sqrt{R} \tag 3.17 \\
\sum_{0<r\leq R}{r d(r)\over \phi(r)} & \ll R \log 2R . \tag
3.18 \endalign $$
\endproclaim
\demo{Proof} To prove (3.16), note
$$ \align
\sum_{r \leq R}{\mu^2 (r)\sigma(r)\over \phi(r)} = &
\sum_{r \leq R}\mu^2 (r) \prod_{p\mid r}(1+{2\over p-1}) =
\sum_{r \leq R}\mu^2 (r) \sum_{d\mid r}{2^{\nu(d)}\over
\phi(d)} \\
& = \sum_{d \leq R}{\mu^{2}(d)2^{\nu(d)}\over \phi(d)}
\sum_{\ell \leq {R\over d}}\mu^{2}(\ell)
\leq R \sum_{d \leq R}{\mu^{2}(d)2^{\nu(d)}\over d\phi(d)}
\\
& \leq R\prod_{p}(1+{2\over p(p-1)}) \ll R.
\endalign $$
The proof of (3.17) is similar, and (3.18) was shown in
\cite{3}.
\enddemo

\head 4. The proof of the Theorem \endhead

In this section we calculate the right-hand sides of (2.6)
and (2.7),
and so obtain our result.
\proclaim{Proposition 1} Assuming GRH, we have
$$\sum_{\sc n \leq N \atop \sc
n \equiv a (\bmod\, q)}\Lambda(n)\lambda_{R}(n) =
{N\log R\over \phi(q)} + {cN\over \phi(q)} + O({N\over
\phi(q)\sqrt{R}})
+ O(N^{{1\over 2}}\log^3 N) + O(R\log N). \tag 4.1$$
\endproclaim
\demo{Proof} Starting from the definition (3.1) and
recalling $L_k(R)$ from
(3.3), we have
$$ \align \sum_{\sc n \leq N \atop \sc
n \equiv a (\bmod\, q)}\Lambda(n)\lambda_{R}(n) & =
\sum_{r \leq R}{\mu^2 (r)\over \phi(r)}\sum_{d \mid
r}d\mu(d)
\sum_{\sc n\leq N \atop {\sc n \equiv a (\bmod\, q) \atop
\sc d\mid
n}}\L(n)\\
& = L_{1}(R)\psi(N;q,a) -
\sum_{r \leq R}{\mu^2 (r)\over \phi(r)}\sum_{p\mid r}p\log p
\sum_{\sc k \geq 1 \atop {\sc p^k \leq N \atop \sc p^k
\equiv a (\bmod\,
q)}}1.
\endalign $$
Here the sum over $k$ is trivially of size
$O(\displaystyle{\log N\over \log p})$, so that by Lemma 3
$$ \align
\sum_{r \leq R}{\mu^2 (r)\over \phi(r)}\sum_{p\mid r}p\log p
\sum_{\sc k \geq 1 \atop {\sc p^k \leq N \atop \sc p^k
\equiv a (\bmod\,
q)}}1
 &\ll \log N \sum_{r \leq R}{\mu^2 (r)\over
\phi(r)}\sum_{p\mid r}p \\
& = \log N \sum_{p \leq R} {p\over \phi(p)}
\sum_{\sc m \leq {R\over p}\atop \sc (m,p)=1}
{\mu^2 (m)\over \phi(m)} \\ &\ll \log N \sum_{p\le R}\log
{R\over p}
\ll R\log N .\endalign$$
By the prime number theorem we obtain (4.1).
In order for the main term to dominate the error terms in
(4.1)
we will require that
$$qR \leq {N\over \log N} ,
\quad q \leq {N^{{1\over 2}}\over \log^3 N} . \tag 4.2 $$
\enddemo
Hence the relevant contribution to (2.6\, b) will be
$$ \align
\left(\int_{2x}^{2x+h}  - \int_{x}^{x+h}\right)&
\sum_{\sc n \leq u \atop \sc
n \equiv a (\bmod\, q)}\Lambda(n)\lambda_{R}(n) \, du \tag
4.3 \\
= & {xh\over \phi(q)}(\log R + c) + O({xh\over
\phi(q)\sqrt{R}})
+  O(x^{{1\over 2}}h\log^3 x) + O(Rh\log x).  \endalign
$$
\proclaim{Proposition 2} Assuming GRH, we have for
$$ 1 \leq q \leq h \leq x , \quad qR \leq x,
\tag 4.4 $$
that
$$ \align \sum_{\sc 0 < k \leq h \atop \sc q \mid k} &
\sum_{\sc x < n \leq 2x+h-k \atop \sc n \equiv a (\bmod\,
q)}
[\Lambda(n)\lambda_{R}(n+k)+ \Lambda(n+k)\lambda_{R}(n)]
f(n,x,h-k) \\
=  & {xh^2 \over \phi^2 (q)}-{xh\over \phi(q)}\log {h\over
q}
+ O({xh\over \phi(q)}(\log\log 3q)^3)
+ O({xh^2 d(q)\over \phi^{2}(q)R}\log{2h\over q}) \\
& +O({x^{{1\over 2}}h^{{3\over 2}}R\log^2 x\over q^{{1\over
2}}})
+O({x^{{1\over 2}}h^{2} R^{{1\over 2}}\log^2 x\over q}) .
\tag 4.5
\endalign$$
\endproclaim
\demo{Proof} We have
$$\sum_{\sc N_{1} < n \leq N_{2} \atop \sc
n \equiv a (\bmod\, q)}\lambda_{R}(n)\Lambda(n+jq) =
\sum_{r \leq R}{\mu^2 (r)\over \phi(r)}\sum_{d \mid
r}d\mu(d)
\sum_{\sc N_1 < n\leq N_2 \atop
{\sc n \equiv a (\bmod\, q) \atop \sc d\mid n}}\L(n+jq).
$$
We may write the innermost sum as
$$
\sum_{\sc N_1+jq < m \leq N_2+jq \atop
{\sc m \equiv a (\bmod\, q) \atop \sc m \equiv jq (\bmod\,
d)}}\L(m).$$
Here $m-jq = \ell d$ for some integer $\ell$, and so $a
\equiv \ell d
(\bmod\, q)$. Since $(a,q) = 1$, we can include only those
$d$'s such that
$(d,q)=1$. Then there is a unique $b, \,\, 0 < b < qd$, such
that $m \equiv
b (\bmod\, qd)$. We know $(m,q)=1$, so that $(m,d)=1$ if and
only if
$(j,d)=1$.
Hence the innermost sum is equal to
$$\align
\psi(N_2 + jq; qd,b) - &  \psi(N_1 + jq; qd,b) \\
& =
{N_2 - N_1 \over \phi(qd)}E_{qd,b} + E(N_2 + jq; qd,b) -
E(N_1 + jq; qd,b),
\endalign $$
where $E_{qd,b} = 1$ if $(qd,b)=1$, and $E_{qd,b} = 0$ if
$(qd,b)>1$.
Thus
$$\align \sum_{\sc N_{1} < n \leq N_{2} \atop \sc
n \equiv a (\bmod\, q)}\lambda_{R}(n) & \Lambda(n+jq)  =
{u-|j|q\over \phi(q)}
\sum_{r \leq R}{\mu^2 (r)\over \phi(r)}
\sum_{\sc d \mid r \atop \sc (d,jq)=1}{d\mu(d)\over \phi(d)}
\tag 4.6 \\
+ &
\sum_{r \leq R}{\mu^2 (r)\over \phi(r)}
\sum_{\sc d \mid r \atop \sc (d,q)=1}d\mu(d)[E(N_2 + jq;
qd,b)
- E(N_1 + jq; qd,b)] ,
\endalign $$
where the first term on the right-hand side is the main
term, its value
settled by Lemma 6, and its contribution to (2.6\, b) will
be
$$ \align
\int_{2x}^{2x+h} & \sum_{0 < |j| \leq {u-2x\over q}}
{u-|j|q\over \phi(q)}[\gs(jq) + O({jqd(jq)\over
R\phi(jq)})]\, du \\
& -\int_{x}^{x+h} \sum_{0 < |j| \leq {u-x\over q}}
{u-|j|q\over \phi(q)}[\gs(jq) + O({jq\, d(jq)\over
R\phi(jq)})]\, du\\
= & {2x\over \phi(q)}\sum_{0<j\leq {h\over q}}(h-jq)\gs(jq)
+
O({xhq d(q)\over R\phi^{2}(q)}\sum_{0<j\leq {h\over q}}{j
d(j)\over
\phi(j)})\\
= & {xh^2\over \phi^2(q)} - {xh\over \phi(q)}\log {h\over q}
+
O({xh\over \phi(q)}(\log\log 3q)^3) +
O({xh^2 d(q)\over \phi^{2}(q)R}\log{2h\over q})
, \tag 4.7 \endalign
$$
by Lemma 7 and (3.18).
For the second term in the right-hand side of
(4.6), if we use (1.4) directly, we will get
the upper bound $Rx^{{1\over 2}}\log^2 x$, by (3.16). This
will
lead to a contribution of $\displaystyle O({x^{{1\over
2}}h^2 R
\over q}\log^2 x)$ in (2.6\, b). Instead, in view of the
averaging
over $j$ in (2.6\, b), we will use Hooley's estimate quoted
as Lemma 8
above.
To do this note that some of the $d$'s may not be coprime to
$b$, but we can
discard them (from the $j$- and $n$-summations) with an
error
$$
\ll \sum_{\sc 0 < |j| \leq {h\over q} \atop \sc
(j,d)>1}\psi(3x;qd,b)
\ll \sum_{\sc 0 < |j| \leq {h\over q}} \sum_{p\mid
d}\sum_{\sc n \leq 3x
\atop
\sc p \mid n}\L(n) \ll {h\over q}\log^2 x \, ,
$$
and this leads to an error of $O(\displaystyle{h^2 R\over q}\log^2 x)$ in
(2.6\, b).
Hence
the contribution to (2.6\, b) from
the second term on the right-hand side of (4.6) is
$$ \align
\ll & {h^2 R\over q}\log^2 x +
\left(\int_{2x}^{2x+h} \!\!- \int_{x}^{x+h}\right)
\sum_{r \leq R}{\mu^2 (r)\over \phi(r)}\!\!
\sum_{\sc d \mid r \atop \sc (d,q)=1}\! d\!\!
\sum_{\sc 0 < |j| \leq {h\over q} \atop \sc
(j,d)=1}\!\!\max_{u\leq 2x+h}
|E(u; qd,b)| \, du \\
\ll &  {h^2 R\over q}\log^2 x + h \sum_{r \leq R}{\mu^2
(r)\over \phi(r)}
\sum_{\sc d \mid r \atop \sc (d,q)=1}d \, ({h\over
q})^{{1\over 2}}
(\sum_{\sc 0 < |j| \leq {h\over q} \atop \sc
(j,d)=1}\max_{u\leq 2x+h}
|E(u; qd,b)|^2 )^{{1\over 2}} \\
\ll &  {h^2 R\over q}\log^2 x + {h^{{3\over 2}}\over
q^{{1\over 2}}}
\sum_{r \leq R}{\mu^2 (r)\over \phi(r)}\!
\sum_{d \mid r}d (1+({h\over qd}))^{{1\over 2}}
(\!\!\sum_{\sc j (\bmod\, d) \atop \sc
(j,d)=1}\!\!\max_{u\leq 3x}
|E(u; qd,b)|^2 )^{{1\over 2}}.\!\!\tag 4.8 \endalign $$
In the last sum as $j$ runs through the reduced residues
modulo $d$, $b$
runs through those elements of the set $\{ a, a+q, \ldots,
a+(d-1)q \}$
which are relatively prime to $d$ (note that $ a\equiv \ell
d (\bmod\, q)$
and $(a,q)=1$ implies $ (a,d)=1$),
and this correspondence is one-to-one.
This is because $m\equiv a (\bmod\, q)$ and $m \equiv jq
(\bmod\, d)$
if and only if $m \equiv n_{1}da + n_{2}jq^2 (\bmod\, qd)$
where $n_i$ satisfy $dn_1 \equiv 1 (\bmod\, q), \; qn_2
\equiv 1 (\bmod\,
d)$,
and we have $n_{1}da + n_{2}jq^2 \equiv a+tq (\bmod\, qd)$
if and only if $j-t \equiv an_{2} (\bmod\, d)$. Hence we may
replace the
$j$-sum in (4.8) by
$$\sum_{\sc t (\bmod\, d) \atop \sc (a+tq,d)=1}\max_{1\leq
u\leq 3x}
|E(u; qd,a+tq)|^2 . \tag 4.9 $$
Although the last sum is over only $\displaystyle{1\over \phi(q)}$ of the
reduced residue
classes modulo $qd$, we shall use Hooley's estimate as is.
One would
want to get a Hooley-type estimate for (4.9) itself, thereby
saving a factor
of $\phi(q)$, but this seems to require some
estimates for certain integrals involving
pairs of $L$-functions. We do not follow this
path now. Recall that
Theorem C, which gives an asymptotic estimate for our
integral
already rests upon such an assumption, (1.10), about
$L$-functions.
By Lemma 8, we take $x\log^4 x$ as upper bound for (4.9) on
the
condition
that $qR \leq x $, and on applying Lemma 9 we obtain that
the expression in (4.8) is
$$
 \ll
{x^{{1\over 2}}h^{{3\over 2}}R\over q^{{1\over 2}}}\log^2 x
+
{x^{{1\over 2}}h^{2}R^{{1\over 2}}\over q}\log^2 x. \tag
4.10
$$
This completes the proof of Proposition 2.
\enddemo
\proclaim{Proposition 3} For $(a,q) = 1$ we have
$$\sum_{\sc n \leq N \atop \sc
n \equiv a (\bmod\, q)}\lambda_{R}^{2}(n) =
{N\over \phi(q)}[\log R + c + O(v(q)) + O(R^{-{1\over
2}+\epsilon})] +
O(R^2).
\tag 4.11 $$
\endproclaim
\proclaim{Proposition 4} For $(a,q) = 1$ and $j\neq 0$ we
have
$$\sum_{\sc n \leq N \atop \sc
n \equiv a (\bmod\, q)}\lambda_{R}(n)\lambda_{R}(n+jq) =
{N\over \phi(q)}\gs(jq) + O({N g(q)\over \phi(q)R}
{j d(j)\over \phi(j)}) + O(R^2).
\tag 4.12 $$
\endproclaim
\demo{Proof} The beginning of the proof of Proposition 3 may
be
incorporated into that of Proposition 4 upon a notational
stipulation
for the case $j=0$. When the positive integer $t$ satisfies
$t\mid j$, if $j=0$ we will understand that $t$ can be any
positive integer;
and we will take $(t,0) = t$.

By definition (3.1),
$$\sum_{\sc n \leq N \atop \sc n \equiv a (\bmod\,
q)}\lambda_{R}(n)
\lambda_{R}(n+jq) =
\sum_{r,\, r' \leq R}{\mu^{2}(r)\mu^{2}(r')\over
\phi(r)\phi(r')}
\sum_{\sc d\mid r \atop \sc e\mid r'}d\mu(d) e\mu(e)
\sum_{\sc n\leq N \atop {\sc n \equiv a (\bmod\, q) \atop
\sc d\mid n ,\,
e\mid n+jq}}1 . \tag 4.13 $$
In the innermost sum, the conditions on $n$, $d$,  and  $e$
imply $(q,de)=1,\,
(d,e)\mid j$, and thus $n$ belongs to a unique residue class
modulo
$[q,d,e]$.
Hence we have
$$ \sum_{\sc n\leq N \atop {\sc n \equiv a (\bmod\, q) \atop
\sc d\mid n ,\,
e\mid n+jq}}1 = {N\over [q,d,e]}+O(1). \tag 4.14 $$
The contribution of the $O(1)$-term in (4.14) to (4.13) is
$$ \ll \sum_{r,\, r' \leq R}{\mu^{2}(r)\mu^{2}(r')\over
\phi(r)\phi(r')}
\sum_{\sc d\mid r \atop \sc e\mid r'}d e =
(\sum_{r\leq R}{\mu^{2}(r)\sigma(r)\over \phi(r)})^2 \ll R^2
\tag 4.15 $$
by (3.16), and this is where the $O(R^2)$-term in (4.11) and
(4.12)
comes from. Hence
$$ \sum_{\sc n \leq N \atop \sc n \equiv a (\bmod\,
q)}\lambda_{R}(n)
\lambda_{R}(n+jq) =
{N\over q}\sum_{r,\, r' \leq R}{\mu^{2}(r)\mu^{2}(r')\over
\phi(r)\phi(r')}
\sum_{\sc d\mid {r\over (r,q)} \atop {\sc e\mid {r'\over
(r',q)} \atop \sc
(d,e)\mid j}}\mu(d)\mu(e)(d,e) \,\, + O(R^2). \tag 4.16 $$
Let $(d,e)= \de, \, d=d'\de, \, e=e'\de$, so that
$(d',e')=1$. The inner
sums
over $d$ and $e$ become
$$
\sum_{\sc \de\mid j \atop \sc \de\mid ({r\over (r,q)},
{r'\over (r',q)})}\de
\sum_{d'\mid {r\over \de(r,q)}}\mu(d')
\sum_{\sc e'\mid {r'\over \de(r',q)} \atop \sc
(e',d')=1}\mu(e'). \tag
4.17$$
Here the innermost sum is
$$\sum_{\sc e'\mid {r'\over \de(r',q)} \atop \sc
(e',d')=1}\mu(e') =
\prod_{\sc p\mid {r'\over \de(r',q)} \atop \sc p \nmid
d'}(1+\mu(p)) =
\left\{ \eqalign{
&1 \quad \hbox{if
${r'\over \de(r',q)} \mid d'$ ;}\cr
&0  \quad\hbox{otherwise.}\cr}\right. \tag 4.18 $$
Next the sum over $d'$ becomes
$$
\sum_{\sc d'\mid {r\over \de(r,q)} \atop \sc {r'\over
\de(r',q)}\mid
d'}\mu(d')
=\left\{ \eqalign{
\mu({r\over \de(r,q)}) &\quad \hbox{if
${r'\over (r',q)} = {r\over (r,q)} $ ;}\cr
0 \qquad &\quad\hbox{otherwise,}\cr}\right. \tag 4.19 $$
so the main term of (4.16) is
$${N\over q}\sum_{\sc r,\, r' \leq R \atop {r\over (r,q)} =
{r'\over
(r',q)}}
{\mu^{2}(r)\mu^{2}(r')\over \phi(r)\phi(r')}\mu({r\over
(r,q)})
\sum_{\de\mid ({r\over (r,q)} , j)}\de\mu(\de) . \tag 4.20
$$
Since
$$\sum_{\de\mid ({r\over (r,q)} , j)}\de\mu(\de) =
\mu(({r\over (r,q)}, j))\phi(({r\over (r,q)}, j))\, , \tag
4.21 $$
the main term is
$${N\over q}\sum_{\sc r,\, r' \leq R \atop {r\over (r,q)} =
{r'\over
(r',q)}}
{\mu^{2}(r)\mu^{2}(r')\over \phi(r)\phi(r')}\mu({r\over
(r,q)})
\mu(({r\over (r,q)}, j))\phi(({r\over (r,q)}, j)). \tag 4.22
$$
Writing $(r,q)=\ell, \, (r',q)=m,\, r=\ell s,\, r'=ms$ where
$(s,q)=1$,
(4.22) takes the form
$$ {N\over q}\sum_{\ell\mid q}{\mu^{2}(\ell)\over
\phi(\ell)}
\sum_{m\mid q}{\mu^{2}(m)\over \phi(m)}
\sum_{\sc s\leq \min ({R\over \ell},{R\over m}) \atop \sc
(s,q)=1}
{\mu(s)\over \phi^{2}(s)}\mu((s,j))\phi((s,j)). \tag 4.23 $$

{\it The $j=0$ case:} We rewrite (4.23) as
$$ {N\over q}\sum_{\ell\mid q}{\mu^{2}(\ell)\over
\phi(\ell)}
\sum_{m\mid q}{\mu^{2}(m)\over \phi(m)}
\sum_{\sc s\leq \min ({R\over \ell},{R\over m}) \atop \sc
(s,q)=1}
{\mu^{2}(s)\over \phi(s)}. \tag 4.24 $$
It is convenient to regard (4.24) as
$$ {N\over q}\sum_{\ell\mid q}{\mu^{2}(\ell)\over
\phi(\ell)}
\sum_{m\mid q}{\mu^{2}(m)\over \phi(m)}
\sum_{\sc s\leq {R\over \ell} \atop \sc
(s,q)=1}{\mu^{2}(s)\over \phi(s)}
-{N\over q}\sum_{\ell\mid q}{\mu^{2}(\ell)\over \phi(\ell)}
\sum_{\sc m\mid q \atop \sc m > \ell}{\mu^{2}(m)\over
\phi(m)}
\sum_{\sc {R\over m} < s \leq {R\over \ell} \atop \sc
(s,q)=1}
{\mu^{2}(s)\over \phi(s)} . \tag 4.25  $$
For the first term of (4.25), we observe that
$$\sum_{m\mid q}{\mu^{2}(m)\over \phi(m)} = {q\over
\phi(q)}, \tag 4.26$$
and by Lemma 3 and Lemma 5 we have
$$ \align {N\over q}& \sum_{\ell\mid q}{\mu^{2}(\ell)\over
\phi(\ell)}
\sum_{m\mid q}{\mu^{2}(m)\over \phi(m)}
\sum_{\sc s\leq {R\over \ell} \atop \sc
(s,q)=1}{\mu^{2}(s)\over \phi(s)}\\
& = {N\over \phi(q)}(\log R + c + v(q)) -
{N\over q}\sum_{\ell\mid q}{\mu^2 (\ell)\over
\phi(\ell)}\log \ell
+ O({Nw(q)\over \phi(q)\sqrt{R}}
\sum_{\ell\mid q}{\mu^2 (\ell)\sqrt{\ell}\over \phi(\ell)})
\\
& = {N\over \phi(q)}(\log R + c) + O({N\over
\phi(q)\sqrt{R}}\prod_{p\mid q}
(1+{1\over \sqrt{p}})(1+{\sqrt{p}\over p-1})).
\tag 4.27 \endalign $$
 The last product has logarithm
$$
\sum_{p\mid q}\log (1+{1\over \sqrt{p}}) + \log
(1+{\sqrt{p}\over p-1}) \,\,
\leq
2\sum_{p\mid q}{1\over \sqrt{p}} + O(1) \ll {\sqrt{\log q}\over \log\log
3q},
\tag 4.28 $$
by (3.6). Hence the first term of (4.25) is
$$ {N\over \phi(q)}(\log R + c) + O({N\over
\phi(q)R^{{1\over 2}-\epsilon}})
\tag 4.29 $$
for any arbitrarily small and fixed $\epsilon > 0$. Using
$$
\sum_{\sc {R\over m} < s \leq {R\over \ell} \atop \sc
(s,q)=1}
{\mu^{2}(s)\over \phi(s)} = {\phi(q)\over q}\log {m\over
\ell} +
O(w(q)\sqrt{{m\over R}}), \tag 4.30 $$
which is implied by Lemma 3, the second term of (4.25) is
expressed as
$$ \align
{N\phi(q)\over q^{2}}\sum_{\ell\mid q}{\mu^{2}(\ell)\over
\phi(\ell)} &
\sum_{\sc m\mid q \atop \sc m > \ell}{\mu^{2}(m)\over
\phi(m)}\log m
- {N\phi(q)\over q^{2}}\sum_{\ell\mid q}{\mu^{2}(\ell)\over
\phi(\ell)}\log \ell
\sum_{\sc m\mid q \atop \sc m > \ell}{\mu^{2}(m)\over
\phi(m)}\\
&  + O({Nw(q)\over
q\sqrt{R}}\sum_{\ell\mid q}{\mu^{2}(\ell)\over \phi(\ell)}
\sum_{\sc m\mid q \atop \sc m > \ell}{\mu^{2}(m)\over
\phi(m)}\sqrt{m} \,).
\tag 4.31 \endalign $$
Each term of (4.31) is majorized by deleting the restriction
$m > \ell$.
Then, by (4.26) and (3.8),
the first two terms are each $\ll \displaystyle{Nv(q)\over
\phi(q)}$,
and the error term is the same as that of (4.27). Hence we
have shown
(4.11).

{\it The $j \neq 0$ case:} Since
$$
\sum_{\sc s=1 \atop \sc (s,q)=1}^{\infty}
{\mu(s)\over \phi^{2}(s)}\mu((s,j))\phi((s,j)) = \prod_{\sc
p\mid j \atop
\sc
p\nmid q}(1+{1\over p-1}) \prod_{p\nmid jq}(1-{1\over
(p-1)^2})
= \gs(jq){\phi(q)\over q}, \tag 4.32 $$
by (3.10) and (3.11), (4.23) is
$$
{N\over \phi(q)}\gs(jq) + O({N\over q}\sum_{\ell\mid q}
{\mu^{2}(\ell)\over \phi(\ell)}\sum_{m\mid
q}{\mu^{2}(m)\over \phi(m)}
\sum_{s > \min ({R\over \ell},{R\over m}) }
{\mu^{2}(s)\mu^{2}((s,j))\phi((s,j))\over \phi^{2}(s)}).
\tag 4.33 $$
The sum over $s$ was encountered before in (3.13) and
majorized
as in (3.9), so the $O$-term in (4.33) is
$$ \align
& \ll {N\over qR}{j\, d(j)\over \phi(j)}
\sum_{\ell\mid q}{\mu^{2}(\ell)\over \phi(\ell)}
\sum_{m\mid q}{\mu^{2}(m)\over \phi(m)}\max (\ell,m) \\
& \ll {N g(q)\over \phi(q) R}{j\, d(j)\over \phi(j)}. \tag
4.34
\endalign $$
This completes the proof of Proposition 4.
\enddemo
We now return to the proof of the theorem. The relevant
contributions to (2.7) are
$$
\sum_{\sc x < n \leq 2x+h \atop \sc n \equiv a (\bmod\,
q)}\lambda_{R}^{2}(n)
f(n,x,h) = {xh\over \phi(q)}(\log R + c + O(v(q))+
O(R^{-{1\over
2}+\epsilon}))
+ O(hR^2) , \tag 4.35 $$
and
$$ \align
& \sum_{\sc 0 < k \leq h \atop \sc q \mid k}
\sum_{\sc x < n \leq 2x+h-k \atop \sc n \equiv a (\bmod\,
q)}
\lambda_{R}(n)\lambda_{R}(n+k) f(n,x,h-k) \tag 4.36 \\
& = \!
{x\over \phi(q)}\sum_{0<j \leq {h\over q}}\gs(jq) +
O({xh\over R\phi(q)}g(q)\sum_{j \leq {h\over q}}{j d(j)\over
\phi(j)}) \\
& = \!{1\over 2}{xh^{2}\over \phi^{2}(q)} -
{1\over 2}{xh\over \phi(q)}\log {h\over q} +
O({xh\over \phi(q)}(\log\log 3q)^3)\!
+O({xh^{2}g(q)\over Rq\phi(q)}\log {2h\over
q})\!+\!O({h^{2}R^{2}\over q}) .
\endalign $$
Now we put together equations  (2.2), (2.5),
(2.6), (2.7), (4.3), (4.5), (4.35), (4.36), subject to
(1.3), (4.2)
and $qR^{2} \leq x$, to obtain
$$\align
I(x,h,q,a) \geq & {xh\over \phi(q)}\log {Rq\over h} +
O({xh\over \phi(q)}(\log\log 3q)^3) +
O({xh^{2}\over R\phi(q)}({d(q)\over \phi(q)}+{g(q)\over q}))
\\
& \,\,\, +O({x^{{1\over 2}}h^{{3\over 2}}R\log^2 x\over
q^{{1\over 2}}})
+O({x^{{1\over 2}}h^{2} R^{{1\over 2}}\log^2 x\over
\phi(q)}) +
O({h^{2}R^{2}\over q}) .  \tag 4.37 \endalign
$$
Recall that $d(q)$ and $g(q)$ are both $\ll q^{\epsilon}$.
Here we pick
$$
R= ({x\over hq\log^4 x})^{{1\over 2}}. \tag 4.38 $$
This choice of $R$ makes all the error terms
$o(\displaystyle{xh\over q})$ provided that
$ h \leq (xq)^{{1\over 3}-\epsilon}$ .

\enddocument